
\magnification=\magstep1
\hsize16truecm
\vsize23.5truecm
\topskip=1truecm
\raggedbottom
\abovedisplayskip=3mm
\belowdisplayskip=3mm
\abovedisplayshortskip=0mm
\belowdisplayshortskip=2mm
\normalbaselineskip=12pt
\normalbaselines
\input mssymb.tex
\def\F{\Bbb F}

\def\P{\Bbb P}

\vskip 2.5pc
\noindent
\vskip 2.0pc
\font\titlefont= cmcsc10 at 12pt
\centerline {\titlefont  Quadratic forms, generalized Hamming weights of
codes} \centerline {\titlefont  and curves with many points.}


\noindent
\vskip 2pc
\font\authorfont= cmcsc10 at 10pt
\centerline
{\authorfont Gerard van der Geer and Marcel van der Vlugt}
\bigskip
\def \a {\alpha}
\bigskip
\noindent
{\bf Introduction.}
\footnote{alg-geom/9412011}
\smallskip
\noindent
In the interaction between algebraic geometry and
classical coding theory  the correspondence between words in trace
codes ${\cal C}$ and  Artin-Schreier curves over finite fields is a
central theme.  This correspondence can be extended to a relation
between subcodes  of ${\cal C}$ and fibre products of Artin-Schreier
curves. As a  consequence weights of words and subcodes are related to
numbers of  rational points on curves, where low weight words and
subcodes yield curves with many rational points. Subcodes in
which all non-zero words have minimum weight lead to curves with a large
number of points compared to known upper bounds.

The most general upper bound for the number of rational
points $\# C (\F_q)$  on a smooth projective
curve $C$ of genus $g$ over a finite field with $q$ elements $\F_q$, is
$$
\# C (\F_q) \leq q + 1 + g [2 \sqrt {q}] \eqno (1)
$$
which is Serre's improved form of the Hasse-Weil upper bound  (see
[Se]). Here $[x]$ is the greatest integer function.

Trace codes which are subcodes of (binary) second-order Reed-Muller codes have
proved to be successful in the construction of
curves with many points. The words in these trace codes correspond not
only to algebraic curves but also to {\sl quadratic forms} over $\F_2$.
Using the well known theory of quadratic forms we shall develop a simple
tool to construct subcodes in such trace codes consisting of minimum
weight words, by which we obtain curves carrying many rational
points.

We shall compare our curves with the curves in the table of
Wirtz [W]. This table lists for many small fields $\F_q$ and genera $g$
an interval [A,B] where $A$ means that there exists a curve $C$ of genus
$g$ over $\F_q$ with $A$ rational points, while $B$ is an upper bound.

Although we restrict to binary codes, a slightly modified
construction can be applied to trace codes over $\F_p$ with odd prime
$p$.

The paper is organized as follows. In Section 1 we recall
elementary facts on fibre products of curves and generalized Hamming
weights of codes and we introduce the codes we shall consider. In
Section 2 the relation between the codes and quadratic forms is
explained. Then we describe the construction of low weight subcodes in
case $q = 2^m$ with $m$ odd in Section 3 and in Section 4 we give
applications of the construction. In the last Section we discuss the
case $q = 2^m$ with $m$ even.

\medskip
\noindent
{\bf \S1. Fibre products of curves and generalized Hamming weights.}
\smallskip
In a series of papers [GV-2, 3, 4] we showed that there is a
relation between generalized Hamming weights of trace codes and the
number of rational points on fibre products of curves associated to
these trace codes. We used  that relation to construct curves with many
rational points. In this section we recall some of the ingredients.

Let $\F_q$ be a finite field of cardinality $q = 2^m$. For $0 < h \leq
[m/2]$ we consider the finite dimensional $\F_q$-vector space of
$2$-{\sl linearized polynomials}
$$
R_h = \{R = \sum_{i = 0}^{h} a_i x^{2^i}\colon  a_i \in \F_q \}.
$$

The vector space  $R_h$ defines a binary linear code ${\cal C}_h$ by
applying the trace map ${\rm Tr}$ from $\F_q$ to $\F_2$:
$$
{\cal C}_h = \{ c_R = ({\rm Tr} (x R(x))_{x \in
\F_{q}^{*}} : R \in R_h\}.
$$
The code ${\cal C}_h$ is a subcode of the punctured binary second-order
Reed-Muller code and for $h = 1, 2$ the codes ${\cal C}_h$ are the duals
of the $2$- and $3$-error correcting BCH codes of length $q - 1$.

To a non-zero word $c_R$ of ${\cal C}_h$ we associate a non-singular
projective  curve $C_R$ with affine equation.
$$
y^2 + y = x R (x).
$$
It is an Artin-Schreier cover of $\P^1$. If $a_h \not = 0 $ the curve
$C_R$ has genus $2^{h - 1}$ and the weight $w (c_R)$ of the word $c_R$ is
related to the number of $\F_q$-rational points on $C_R$, by the
equation $$ w (c_R) = q -  (\# C_R (\F_q) - 1 ) / 2.
$$

If ${\cal D} \subset {\cal C}_h$ is a $r$-dimensional subcode we associate a
curve
$C^{({\cal D})}$ to ${\cal D}$ as follows. Choose a $\F_2$-basis
$c_{R_1},..., c_{R_r}$ of ${\cal D}$ and set
$$
C^{({\cal D})} = ~{\rm Normalization~of}~ C_{R_1} \times_{\P^1} \times ...
\times_{\P^1} C_{R_r},
$$
where the fibre product is taken with respect to the canonical projections
$\varphi_{R_i} : C_{R_i} \to \P^{1}$. Up to isomorphism this does not
depend on the chosen basis.
\medskip
In [G-V2] we proved the relation
$$
t_{\rm Frob} C^{({\cal D})} = \sum_{c_{R_i} \in {\cal D} - \{0\}}
t_{\rm Frob} (C_{R_i}), \eqno (2)
$$
between the traces of Frobenius of the curves involved. Here  the trace
of Frobenius for a curve $C$ over $\F_q$ is defined as
$t_{\rm Frob} = q + 1- \#C (\F_q)$. Furthermore, we have the relation for
the genera $g$ of the curves involved:
$$
g (C^{({\cal D})}) = \sum_{c_{R_i}
\in {\cal D} - \{0\} } g (C_{R_i}). \eqno (3)
$$

An important parameter of a subcode ${\cal D}$ of a code ${\cal C}$ is
its {\sl weight} $w({\cal D})$, by which we mean the number
of coordinate places for which at least one word of ${\cal D}$ has a non-zero
coordinate. In the binary case we have the
relation
$$
w ({\cal D})	= {1 \over 2^{r-1}} \sum_{d \in {\cal D}} w (d).
$$

The {\sl r-th generalized Hamming weight} $d_r ({\cal C})$ of ${\cal C}$ is
defined by
$$
d_r ({\cal C}) = {\rm min}  \{ w ({\cal D}): {\cal D} ~~{\rm is~
a}~ r {\rm -dimensional~ subcode ~of}~~{\cal C} \}.
$$

If there exists a $r$-dimensional subcode in which all non-zero words
have minimum weight, which we shall call a {\sl minimum weight subcode},
then $$
d_r ({\cal C}) = (2^r - 1) d_1 ({\cal C}) / 2^{r-1}. \eqno (4)
$$

In our case where ${\cal C} = {\cal C}_h$ we recall the following proposition
from [G-V2].

\medskip
\proclaim
(1.1) Proposition.  The weight $w ({\cal D})$ of the $r$-dimensional
subcode ${\cal D} \subset {\cal C}_h$ satisfies
$$w ({\cal D}) = q - (\# C^{({\cal D})} (\F_q) - 1)/2^r.
\hfill\qquad \square$$ \par

We immediately see that subcodes of low weight correspond to fibre products
with many points.
\bigskip
\noindent
{\bf \S2. Quadratic forms and codes.}
\smallskip
To a non-zero word $c_R \in {\cal C}_h$ (or to $R \in R_h - \{0\}$) we
can associate not only  a curve
but  also a {\sl quadratic form}. Indeed, the expression
$$
Q_R (x) = {\rm Tr} (x R (x)),
$$
defines on the $\F_2$-vector space $\F_{q = 2^m}$  a quadratic form
over $\F_2$ in $m$ variables. In the sequel we  denote $Q_R (x)$
simply by $Q(x)$. The associated symmetric bilinear form
$$
B (x, y) = {\rm Tr} (x R (y) + y R (x)),
$$
which is also symplectic, has {\sl radical}
$$
W = \{x \in \F_q: B (x, y) = 0 \quad {\rm for~all} \quad  y \in \F_q\}.
$$
This $\F_2$-vector space $W$ has dimension $w$ with $m - w \equiv 0 ({\rm
mod}~2)$ since $B$ is symplectic. If $a_h \not= 0$ then $0 \leq w \leq
2h$ (see [G-V1]). Quadratic forms $Q$ and $Q^\prime$ in $\F_2 [x_1,...,
x_m]$ which can be obtained from each other by a coordinate
transformation are called {\sl equivalent quadratic forms}, denoted by
$Q \sim Q^\prime$. We define the {\sl rank} of $Q$ by
$$
{\rm rk}(Q) = {\rm min}_{Q^\prime \sim Q} ~ \{ {\rm number
{}~of~variables~actually~occuring~in}~ Q^\prime \}.
$$

In characteristic $2$ the quadratic form $Q(x)$ is not necessarily zero
on the radical $W$. Therefore we introduce $$W_0 = \{ x \in W: Q (x)
= 0 \}$$ and we recall from [G-V1] that ${\rm dim} (W_0) = {\rm
dim} (W)$ or ${\rm dim} (W_0) = {\rm dim} (W) - 1$. The theory of
quadratic forms in characteristic 2 (see [A]) yields
\medskip
\proclaim  (2.1) Proposition. The quadratic form $Q (x) = {\rm Tr} (x R
(x))$, with $R \in R_h - \{0\}$, has rank $m - w$ or $m - w + 1$
according to $W_0 = W$ or $W_0 \not= W$. \hfill\qquad $\square$ \par
\medskip
Over $\F_2$  the classification of quadratic forms $Q (x)$
 is strongly related to the number of zeros $N (Q) = \# \{x \in
\F_q: Q (x) = 0 \}$ (see [L-N, Ch. 6]). Indeed, we have:
\medskip
\noindent
\proclaim  (2.2) Proposition. If $Q (x)$ has rank $m - w + 1$
(odd) then
$$
\eqalign{
Q (x) &\sim X_1 X_2 + ...+ X_{m - w - 1} X_{m -w} + X_{m - w
+ 1}^{2} \quad {\rm and} \cr
N (Q) &= 2^{w - 1} \cdot 2^{m - w} =q/2.
}
$$
If $Q (x)$ has rank $m - w $ (even) then either
$$
\eqalign{
Q (x) &\sim X_1 X_2 + ... + X_{m - w - 1} X_{m - w} {\rm ~~and}\cr
N(Q) &= 2^{w} (2^{m - w - 1} + 2 ^{(m - w - 2)/2}) =
(q + \sqrt{q 2^w})/ 2, \cr
}
$$
or
$$
\eqalign{
Q (x) &\sim X_1 X_2 + ... + X_{m-w-1} X_{m-w} + X_{m-w-1}^2 +
X_{m - w}^2
{\rm ~~and }\cr
N (Q) &= 2^w (2^{m - w - 1} - 2^{(m -w-2)/2} ) =
(q - \sqrt {q 2^w} ) /2.\cr
}
$$
\hfill\qquad $\square$ \par
\smallskip

We have simple relations for the weight of a word $c_R \in {\cal C}_h$ and the
number of rational points on the corresponding
curve:
$$
w (c_R) = q - N ({\rm Tr}(x R(x) ) ) \quad {\rm and} \quad \# C_R (\F_q) = 2N
({\rm Tr} (x R (x))) + 1. \eqno (5)
$$

\noindent
{\bf (2.3) Remark.} Words  in ${\cal C}_h$ with $w=2h$ or
$w=2h-1$  correspond to $R$ of degree $2^h$, which means that the
corresponding curves have genus $2^{h-1}$.

\bigskip
\noindent
{\bf \S3. The construction of low weight subcodes of ${\cal C}_h$ for
odd $m$}.
\smallskip
\noindent
In this section we take $q = 2^m$ with $m$ odd, $m
\geq 3$. To find minimum weight subcodes in ${\cal C}_h$ with $0 < h
\leq (m - 1) /2$  we exploit the following simple observation.
\medskip
\proclaim
(3.1) Proposition. For $0 \leq w \leq 2 h - 1$ with $m - w \equiv 0
({\rm mod}~2)$ the expression
$$
\sum_{i = 1}^{(m - w)/2}{\rm Tr} (a_i x) {\rm Tr} (b_i x) \eqno (6)
$$
with $\{a_i,b_i: i = 1,..., (m - w) /2\} \subset
\F_q$ independent over $\F_2$, is  a quadratic form which is equivalent
to $X_1 X_2 + ... + X_{m - w-1} X_{m-w}$ and which has  $(q + \sqrt {q
2^w})/2$ zeros in $\F_q$. \par
\noindent
{\sl Proof.} Since ${\rm Tr}(ax) \in \F_2 [x_1,..., x_m]$ is a linear
form,   (6) is indeed a quadratic form. The $\F_2$-independence of
$\{a_i, b_i\}$ implies that the matrix of the transformation $X_{2i - 1}
= {\rm Tr}(a_i x)$ and $X_{2i} = {\rm Tr} (b_i x)$ for $i = 1,..., (m - w) / 2$
has
rank $m - w$. This matrix can be completed to a coordinate
transformation in the $m$-dimensional $\F_2$-vector space $\F_q$  which
yields the equivalence. The number of zeros follows from Proposition
(2.2). \hfill\qquad $\square$
\medskip
We shall denote the quadratic
form (6) by $Q (a_1,..., a_{(m - w)/ 2}, b_1,..., b_{(m - w)/2})$ or
simply by $Q (a, b)$. The form $Q (a, b)$ evaluated on $\F_{q}^*$
induces a word of length $q - 1$ and often we shall identify the word
and the form. If the $\F_2$-rank of the set
$$
\{a_1,..., a_{(m - w)/2},
b_1,..., b_{(m - w)/2}\},
$$
which we shall also shortly denote by $\{a, b\}$, is $m - w$ then $Q (a,
b)$ has weight $(q - \sqrt {q 2^w})/2$.

First we derive a criterion for words induced by $Q (a,
b)$ to be in ${\cal C}_h$.
\medskip
\noindent
\proclaim (3.2) Proposition. If the
elements $a_i, b_i \in \F_q$ with $1 \leq i \leq (m - w)/2$ satisfy the
system of equations
$$
\sum_{i = 1}^{(m - w)/2}  (a_{i}^{2^j} b_i +
a_i b_{i}^{2^j}) = 0 \eqno (7)
$$
for $j = h + 1, ..., (m - 1)/2$,  then the word of length $q - 1$ induced
by evaluating $\sum_{i = 1}^{(m - w)/2} {\rm Tr} (a_i x) {\rm Tr} (b_i
x)$ on $\F_{q}^*$ is a codeword in ${\cal C}_h$.
	\par
\noindent
{\sl Proof.} Observe that for $a, b \in \F_q$ we
have
$$
{\rm Tr} (ax) {\rm Tr}(bx) =  {\rm Tr} ({\rm Tr} (ax) bx) = {\rm
Tr}  (\sum_{j = 0}^{m - 1} a^{2^j} bx^{2^j +1}	).
$$
The trace map on $\F_q$ satisfies
$$
{\rm Tr}
 (a^{2^j} b x^{2^j +1})= {\rm Tr}  (a b^{2^{m-j}} x^{2^{m - j}+1})
$$
which
implies that
$$
{\rm Tr}(ax){\rm Tr}(bx) = {\rm Tr} (\sum_{j = 0}^{(m -
1)/2} (a^{2^j} b + a b^{2^j}) x^{2^j + 1}).
$$
The argument of ${\rm Tr}$ must be of the form $x R (x)$ with $R \in
R_h$, hence we do not want terms $x^{2^j + 1}$ with $j > h$ to appear in
the argument of ${\rm Tr}$. This condition yields the system of
equations (7). \hfill\qquad $\square$
\medskip
\proclaim
(3.3) Corollary.  For fixed $a_1,..., a_{(m - w)/2}$ the
words $Q (a, b)$ induced by the solutions $(b_1,..., b_{(m - w)/2})$ of
(7) form a subcode of ${\cal C}_h$.
\par
\noindent
{\sl Proof.} From the linearized character of (7) in the $b_i$ we
immediately see that the set of solutions $(b_1,..., b_{(m-w)/2}) \in
\F_{q}^{(m - w)/2}$ is a $\F_2$-subspace of $\F_{q}^{(m - w)/2}$.
\hfill \qquad $\square$
\medskip

For the sequel we fix a
$\F_2$-independent subset $\{a_i : 1 \leq i \leq (m - w) /2 \}$ of
$\F_q$ and we may assume $a_1 = 1$, which is not a restriction because
the equations in (7) are homogeneous. Furthermore, we assume that the
system of equations is not empty. Since we are especially interested in
subcodes consisting of minimum weight words we take $w = 2 h - 1$ and we
put $(m - w) / 2 = (m - 2 h + 1) / 2 = M$. In this situation we have $M
- 1$ equations with $M$ unknowns $b_i$.
\medskip
\noindent
\proclaim (3.4)
Proposition. For fixed $a_1, a_2,..., a_M$ in $\F_q$ the system of
equations (7) has at least $q$ solutions $(b_1, b_2,..., b_M)$ in
$\F_{q}^M$.
	\par
\noindent
{\sl Proof.} If we take an equation
$b_{1}^{2^j} +b_1 = \sum_{i = 2}^{M} a_i b_{i}^{2^j} + a_{i}^{2^j}
b_{i}$ from (7) and if we apply the transformation $b_1 = b_{1}^\prime +
\sum_{i = 2}^{M} \root{2^j} \of {a_i} b_i$ then this equation becomes
$$
b_{1}^{\prime^{2^j}}+ b_{1}^\prime = \sum_{i = 2}^{M}
\bigl( a_{i}^{2^j} + \root{2^j} \of {a_i} \bigr) b_i. \eqno (8)
$$
Note that (8) has $q^{M-1}$ solutions if at least one of the
coefficients of the $b_i$ on the right side is non-zero. Otherwise the
number of solutions depends on the number of solutions of $X^{2^j} + X =
0$  in $\F_q$ and is at least $2 q^{M-1}$. So the dimension of the
$\F_2$-vector space of solutions of an equation from (7), with $m M$
unknowns over $\F_2$, is at least $m (M - 1)$. Since we have $M - 1$
equations our assertion follows. \hfill\qquad $\square$
\medskip
The actual number of solutions of our system depends on the choice of the
$a_i$ as we shall see in examples in the next section. To solve the
system (7) we repeatedly carry out transformations as described in the
proof of Proposition (3.4) and then finally we arrive at an equation of
the form
$$
S (y) = R (x),
$$
where $S$ and $R$ are $\F_2$-linearized
polynomials with coefficients in $\F_q$. This equation can be solved by
lowering the degree step by step. If degree $(S) \leq$ degree $(R)$ we
start with the substitution $y = y^\prime + f (x)$ with a suitable
$2$-linearized polynomial $f$. As long as the resulting polynomials are
not zero we continue the descent until we arrive at a linear equation in
two unknowns which has $q$ solutions.
\smallskip
To obtain words of  minimum weight in ${\cal C}_h$ we must have
solutions  $(b_1,..., b_M)$ of (7) where ${\rm rk}_{\F_2}( \{ a, b\}
)= 2M$.
\medskip
\noindent
\proclaim (3.5) Proposition. Let  $(b_1,..., b_M)$ be a solution of (7).
Then we have:
$$
{\rm rk}_{\F_2} (\{a_1,..., a_M, b_1 \})= M +1
\quad {~\sl implies~~} {\rm rk}_{\F_2}(\{a_1,..., a_M, b_1,..., b_M\})=
2M. $$
	\par
\noindent
{\sl Proof.} Denote $\sum_{i = 1}^{M} {\rm Tr} (a_i x) {\rm Tr} (b_i x)
\in {\cal C}_h$ by $\sum_{ i = 1}^M X_{2i-1} X_{2i}$. Assume that $\{a,
b\}$ has $\F_2$-rank $< 2M$ with $b_M \in ~< a_1,.., a_M, b_1,...,
b_{M-1} > ~$, the $\F_2$-subspace of $\F_q$ generated by the elements
between the pointed brackets. This implies $X_{2M} = \sum_{i = 1}^{2M -
1} \alpha_{i} X_i$ with $\alpha_i \in \F_2$ not all zero. Then
$$
\sum_{i = 1}^{M} X_{2i-1}  X_{2i}=
\sum_{i = 1}^{M - 1} (X_{2i-1} +\alpha_{2i} X_{2M - 1} ) (X_{2i} +
\alpha_{2i- 1} X_{2M-1}) +
$$
$$\qquad \qquad + ( \a_{2M-1} + \sum_{i = 1}^{M - 1} \alpha_{2i
- 1} \alpha_{2i} ) X_{2M-1}^2.
$$
Note that $X_{2M - 1}^{2} = ({\rm Tr} (a_M x))^2 = ({\rm Tr} (a_M x)^2)
\in {\cal C}_h$ and so
$$
\sum_{i = 1}^{M - 1} (X_{2i - 1} + \alpha_{2i} X_{2M - 1}) (X_{2i} +
\alpha_{2i - 1} X_{2M - 1}) = \sum_{i = 1}^{M - 1}{\rm Tr} (\alpha_{i}^{(1)}
x) {\rm Tr}(b_{i}^{(1)} x) \in {\cal C}_h.
$$
Here $a_{i}^{(1)} = a_i + \alpha_{2i} a_M$ and $b_{i}^{(1)} = b_i +
\alpha_{2i -1} a_M$ for $1 \leq i \leq M - 1$ while
$$
{\rm rk}_{\F_2}( \{a_{i}^{(1)}, b_{i}^{(1)}: 1 \leq i \leq M - 1\})\leq
2M - 2.
$$
If ${\rm rk}_{\F_2}( \{a_{i}^{(1)}, b_{i}^{(1)}: 1 \leq i \leq M -
1\} ) = 2 M - 2$ we have found a word in ${\cal C}_h$ of weight $ (q -
\sqrt {q  2^{2 h + 1}})/ 2$ which is smaller than the minimum
weight. On the other hand if the $\F_2$-rank  of $\{a_{i}^{(1)},
b_{i}^{(1)}: 1 \leq i \leq M-1 \} < 2 M - 2$ we repeat the foregoing
procedure. In the end we arrive at a non-zero word $${\rm Tr} ((a_1 +
L_1 (a_2,..., a_m)) x) {\rm Tr} (b_1 + L_2 (a_2 ..., a_m)) x) \in {\cal
C}_h,$$ with $\F_2$-linear combinations $L_1$ and $L_2$ of $a_2,...,
a_m$, which also violates the minimum weight. \hfill\qquad $\square$
\medskip
Let $S$ be the $\F_2$-vector space of solutions $(b_1,..., b_M)$ of (7)
and $V$ the image in $\F_q$ of the projection of $S$ on the first
coordinate $b_1$.
\medskip
\proclaim (3.6) Theorem. If $r = {\rm dim}_{\F_2}(V) - M > 0$ then there
exists a minimum weight subcode of ${\cal C}_h$ of dimension $r$.
	\par
\noindent
{\sl Proof.} Since $r = {\rm dim}_{\F_2} (V) - M > 0$ we can
find $b_{1}^{(1)},..., b_{1}^{(r)}$ in $V$ such that
$$
{\rm rk}_{\F_2}(\{a_1,..., a_M,
b_{1}^{(1)},..., b_{1}^{(r)}\})= M+r.
$$
By Proposition
(3.5) the corresponding solutions $\bigl (b_{1}^{(i)}, b_{2}^{(i)},...,
b_{M}^{(i)}\bigr)$ for $i = 1,..., r$ induce words of minimum weight in
${\cal C}_h$. These words generate a $r$-dimensional subcode of minimum
weight.
	\hfill\qquad $\square$
\medskip
{}From Proposition (3.4) we derive that ${\rm dim}_{\F_2} (S) \geq m$. The
following proposition provides us with a lower bound for ${\rm
dim}_{\F_2}(V)$.
\medskip
\proclaim (3.7) Proposition. Let $(b_1,
b_2,..., b_M) \in S$ then $(b_1, b_{2}^{\prime},..., b_{M}^{\prime}) \in
S$ if and only if $(b_{2}^{\prime},..., b_{M}^{\prime}) = (b_2,..., b_M)
+ (a_2 ,..., a_M)A$, where $A$ is a symmetric $(M - 1) \times (M - 1)$
matrix over $\F_2$. \par
\noindent
{\sl Proof.} The `if-part' follows
immediately by substitution in (7). For the `only if -
part' we observe that
$$
{\rm Tr} (a_i x){\rm Tr} ((b_i + b_{i}^{\prime}) x) \eqno (9)
$$
represents a word in ${\cal C}_h$. Just as in the proof of Proposition
(3.5) this implies that for $i = 2,..., M$ the $(b_i + b_{i}^\prime)$
are $\F_2$-dependent of $\{a_2,..., a_M\}$ or $b_{i} + b_{i}^\prime =
\sum_{j = 2}^{M} \alpha_{ij} a_j$.
\smallskip
By substitution in (9) we obtain
$$
\sum_{i = 2}^{M} \a_{ii} ({\rm Tr} (a_i x))^2 + \sum_{i = 2}^{M-1} {\rm Tr}
(a_i x)
{\rm Tr} ((\sum_{i < j \leq M} \beta_{ij} a_j) x ) \eqno (10)
$$
with $\beta_{ij} = \alpha_{ij} + \alpha_{ji}$. The second term in (10)
also represents a word in ${\cal C}_h$. Again to avoid contradiction
with the minimum weight in ${\cal C}_h$ we must have that $\sum_{i < j
\leq M} \beta_{ij} a_j$ is dependent of $a_2,..., a_{M - 1}$ for $i =
2,..., M$. We then find successively : $\beta_{i, M} = 0$ for $2 \leq i
\leq M - 1, \beta_{i, M - 1} = 0$ for $2 \leq i \leq M - 2$ until
$\beta_{23} = 0$. This means that the $(M - 1) \times (M - 1)$ matrix $A
= (\alpha_{ij})$ is symmetric. \hfill\qquad $\square$
\smallskip
Since we can choose $M (M - 1) / 2$ entries in $A$ independently we
have:
\smallskip
\noindent
\proclaim (3.8) Corollary. The dimension of $V$ over $\F_2$ satisfies
$${\rm dim}_{\F_2} (V) = {\rm dim}_{\F_2} (S) - {M \choose 2} \geq m - {M
\choose 2}. \hfill\qquad \square$$
\par
\smallskip
In the following section we shall illustrate the ideas and results of
this section by some examples.
\bigskip
\noindent
{\bf \S4. Determination of generalized Hamming weights and curves
with many points.}
\smallskip
\noindent
We have

$q = 2^m \quad {\hbox{ \rm  with } }  m \geq 3 \quad {\rm  odd},$

$0 < h \leq (m - 1) / 2,$

$w = 2h - 1\quad {~\rm  and }$

$M = (m - 2 h + 1) / 2.$

\smallskip
\noindent
I. {\sl The case} $h = (m - 1)/2.$
\smallskip
\noindent
Since the system (7) is empty, Propositions (3.1) and (3.2)
imply that  ${\rm Tr}(x){\rm Tr} (bx) \in {\cal C}_h$ and has minimum
weight for $b \in \F_q - \F_2$. By choosing $\F_2$-independent elements
$b^{(1)},..., b^{(m - 1)}$ in $\F_q$ we find
words which generate a $(m - 1)$-dimensional subcode of ${\cal C}_h$
of minimum weight. Combining (4) and Proposition (2.2) we find:
\smallskip
\proclaim
(4.1) Proposition. The generalized Hamming weights of the binary code
${\cal C}_{(m - 1)/2}$ of length $2^m - 1$ satisfy
$$
d_r ({\cal C}_{(m - 1)/2}) = (2^r - 1) \cdot  2^{m - r-1} \quad
{\sl for} \quad 1 \leq r \leq m - 1. \hfill\qquad \square
$$
\par
\smallskip
Applying the fibre product construction from \S1 to these minimum weight
subcodes we get:
\smallskip
\proclaim
\bf (4.2) Corollary. For $1 \leq r \leq m - 1$ there exist curves $C_r$
defined over $\F_q$ of genus $g (C_r) = (2^{r} - 1) 2^{(m - 3)/2}$ and
$\# C_r (\F_q) = 2^m + 1 + (2^r - 1) 2^{m - 1}.$ \par
\noindent
{\sl Proof.} The curves corresponding to the non-zero words in the
$r$-dimensional minimum weight subcode ${\cal D}_r$ have genus $g = 2^{h
- 1} = 2^{(m - 3)/2}$ and according to (5) $t_{\rm Frob} = q - 2N (Tr (x
R (x)) = - \sqrt {q  2^{2h - 1}}$.  \smallskip
{}From (2) and (3) we obtain $g (C^{({\cal D}r)}) = (2^r - 1) 2^{(m - 3)/2}$
and
$$
\#C^{({\cal D}r)} (\F_q) = q + 1 + (2^r - 1) \sqrt {q  2^{2h-1}} =
2^m + 1 + (2^r - 1) 2^{m - 1}. \hfill\qquad \square
$$
This result confirms our earlier result [G-V3, Thm. 5].
\medskip
\noindent
II. {\sl The case } $h = (m - 3) /2$.
\smallskip
We take $m \geq 5$. In this case $M = 2$; fix $\{1,
a_2\}$ with $\F_2$-rank $2$  in $\F_q$. The system (7) consists of the
equation
$$
b_{1}^{2^{(m - 1)/2}} + b_1 = a_2 b_{2}^{2^{(m - 1)/2}} +
a_{2}^{2^{(m - 1)/2}} b_2. \eqno (11)
$$
The transformation $b_1 = b_{1}^{\prime} + \root 2^{(m - 1)/2} \of {a_2}
b_2$ changes (11) into
$$
b_{1}^{\prime{2^{(m - 1)/2}}} + b_{1}^\prime =
\bigl (a_{2}^{2^{(m - 1)/2}} + \root 2^{(m - 1)/2} \of {a_2} \bigr)
b_2. \eqno (12)
$$

The coefficient of $b_2$ in (12) is non-zero because $\{1, a_2\}$ has
$\F_2$-rank $2$. Hence (11) has $q$ solutions
$(b_1, b_2)$ or ${\rm dim}_{\F_2} (S) = m$.
{}From Corollary (3.8) we deduce ${\rm dim}_{\F_2} (V) = m - 1$ and we
find by Theorem (3.6):
\smallskip
\proclaim (4.3) Proposition. i) The generalized Hamming weights of ${\cal
C}_{(m -3)/2}$ of length $2^m - 1$ satisfy
$$
d_r ({\cal C}_{(m - 3)/2})
= (2^r - 1) (2^{m - 1} - 2^{m - 3}) /2^{r-1} \quad {\sl for} \quad 1
\leq r \leq  m - 3.
$$
ii) For $1 \leq r \leq m - 3$ there exist curves
$C_r$ defined over $\F_{2^m}$, with  $m \geq 5$ odd, of genus $g (C_r) =
(2^r - 1) 2^{(m - 5) /2}$ and with $\# C_r (\F_{2^m}) = 2^m + 1 + (2^r -
1) 2^{m - 2}$.
\par
\noindent
{\sl Proof.} The proof is similar to the proof of
Proposition (4.1) and Corollary (4.2).
	\hfill\qquad $\square$
\smallskip
\noindent
{\bf (4.4) Example.} For $m = 7$ we have $h = 2$ and ${\cal C}_2 = BCH
(3)^{\bot}$, the dual of the $3$-error correcting BCH code of length
$127$. The former Proposition yields:
\smallskip

$\bullet$ {\sl The generalized Hamming weights} $d_r, 1 \leq r \leq
4$, {\sl of} $BCH (3)^{\bot}$ {\sl of length} $127$ {\sl are}
$$d_1 = 48,~~d_2 =
72, ~~d_3= 84, ~~d_4 = 90.
$$

$\bullet$ {\sl There exists curves} $C_r$ {\sl
defined over} $\F_{128}$ {\sl with} \medskip
\settabs 5 \columns
\+&$r$ &$g (C_r)$ &$\# C_r (\F_{128})$ &{\rm Wirtz} \cr
\smallskip
\+&$1$ &$2$ &$161$ &$(184-195)$\cr
\smallskip
\+&$2$ &$6$ &$225$ &$(225-261)$\cr
\smallskip
\+&$3$ &$14$ &$353$ &$(289-437)$\cr
\smallskip
\+&$4$ &$30$ &$609$ &$(369-789)$.\cr
\medskip
\medskip
\noindent
III. {\sl  The case } $h = (m - 5) / 2$.
\smallskip
We take $m \geq 7$. In this case $M = 3$; fix $\{1, a_2, a_3\}
{}~\F_2$-independent in $\F_q$. The system of
equations is:
$$b_{1}^{2^{(m - 3)/2}} + b _1 = \sum_{i = 2, 3} \biggl (a_i b_{i}^{2^{(m -
3)/2}} + a_{i}^{2 ^{(m - 3)/2}} b_i \biggr),
\eqno (13)$$

$$b_{1}^{2^{(m - 1)/2}} + b_1 = \sum_{i = 2, 3} \biggl (a_i b_{i}^{2^{(m -
1)/2}} + a_{i}^{2^{(m - 1)/2}} b_i\biggr). \eqno
(14)$$

Apply the linear transformation $b_1 = b_{1}^\prime +
\root s \of {a_2} b_2 + \root s \of {a_3}
b_3$ with $s=2^{(m - 1)/2}$. Equation (14) becomes linear in $b_2$ and
$b_3$ with non-zero coefficients. We express $b_3$ in $b_{1}^\prime$ and
$b_2$. Then (13) becomes an equation in $b_{1}^{\prime}$ and $b_2$ of
the form $S (b_{1}^{\prime})  = R(b_2)$ with $R$ and $S~2$-linearized.
In general $S (b_{1}^\prime)$ has degree $2^{m - 2}$ and $R (b_2)$ has
degree $2^{(m - 3)/2}$ and we find $q$ solutions $(b_1,b_2, b_3)$ (see
the remarks after Proposition (3.4)). Just as in the cases I and II we
find a proposition on generalized Hamming weights and on curves.
\smallskip
The existence of $q$ solutions is guaranteed by Proposition
(3.4). However for special choices of $\{1, a_2, a_3\}$ the number of
solutions is a multiple of $q$. To illustrate this remark we consider
the following example.
\smallskip
\noindent
{\bf (4.5) Example.} For $m = 9$ and $h
= 2$ our system is:
$$
b_{1}^8 +  b_1 = a_2 b_{2}^8 + a_{2}^8 b_2 + a_3
b_{3}^8 + a_{3}^{8} b_3, \eqno (15)
$$
$$
b_{1}^{16} + b_1 = a_2
b_{2}^{16} + a_{2}^{16} b_2 + a_{3} b_{3}^{16} + a_{3}^{16} b_3. \eqno
(16)
$$

Now we choose a special $\F_2$-independent triple: $\{1, a_2, a_3\}
\subset \F_8$. Then we can write (15) as
$$
b_{1}^8 + b_1 = a_{2}^{8} b_{2}^{8} + a_2 b_2 + a_{3}^{8} b_{3}^{8} + a_3
b_3. \eqno (17)
$$

The solutions of (17) are $\{(b_1, b_2, b_3) \in \F_{512}^3: b_1 + a_2 b_2 +
a_3 b_3 = \mu \in \F_8\}$. Take a
value
$\mu \in \F_8$ and substitute $b_1 = \mu + a_2 b_2 + a_3 b_3$ in
(16). We find an equation in $b_2$ and $b_3$ which has $q$ solutions
$(b_2, b_3)$, so the original system has $8q$ solutions and ${\rm
dim}_{\F_2} (S) = 12$. Then Corollary (3.8) and Theorem (3.6) yield:
\smallskip

$\bullet$  {\sl The generalized Hamming weights} $d_r$ {\sl
of} $BCH (3)^\bot$ {\sl of length} $511$ {\sl satisfy}
$$
d_r = (2^r - 1) d_1 /
2^{r -1} = 224 (2^r-1)/2^{r-1} \quad {\sl for} \quad 1 \leq r \leq 6.
$$

$\bullet$ {\sl There exist curves} $C_r$ {\sl defined over}
$\F_{512}$ {\sl with} \medskip
\settabs 5 \columns
{
\+&$r$ &$g({\cal C}_r)$ &$\# C_r (\F_{512})$ &{\sl Upper~bound ~(1)}\cr
\smallskip
\+&$1$ &$2$ &$577$ &$603$\cr
\smallskip
\+&$2$ &$6$ &$705$ &$783$\cr
\smallskip
\+&$3$ &$14$ &$961$ &$1143$\cr
\smallskip
\+&$4$ &$30$ &$1473$ &$1863$.\cr
}
\medskip
\noindent
Combination of the cases also produce good curves as we shall show by an
example.
\smallskip
\noindent
{\bf (4.6) Example.} Take $m = 7$ and fix $\F_2$-independent $\{1, a_2\}$
in $\F_{128}$. By Proposition (4.3) in case II there exists a word of
minimum weight in ${\cal C}_2$, which is $48$, of the form ${\rm Tr} (x)
{\rm Tr} (b_1 x) + {\rm Tr} (a_2 x) {\rm Tr} (b_2 x)$. Now $c_1 = {\rm Tr} (x)
{\rm Tr} (b_1 x)$ and $c_2 = {\rm Tr} (a_2 x) {\rm Tr} (b_2 x)$ are words in
${\cal C}_3$ since in the corresponding ${\rm Tr} (x R (x))$ the degree of $x
R (x)$ is equal to $2^{(m - 1)/2} + 1 = 9$. From Proposition (2.2) it
follows that $c_1$ and $c_2$ have weight $32$. The curves corresponding
to $c_1, c_2$ and $c_1 + c_2$ have genus $4, 4, 2$ respectively and
trace of Frobenius $- 64, -64, -32$ respectively. So the fibre product
curve corresponding to the subcode of ${\cal C}_3$ generated by $c_1$ and
$c_2$ is a curve over $\F_{128}$ of genus $g = 10$ and $129 + 160 = 289$
rational points. The upper bound (1) in this case is $349$.

Applying the same approach to a word
$$
{\rm Tr} (x){\rm Tr} (b_1 x) + {\rm Tr} (a_2 x) {\rm Tr} (b_2 x) + Tr
(a_3 x) Tr (b_3 x)
$$ of minimum weight in ${\cal C}_1$ over $\F_{128}$
one finds a curve of genus $g = 9$ with $241$ rational points (Wirtz:
209-327).
	\bigskip
\noindent
{\bf \S5. The construction of low weight subcodes of
${\cal C}_h$ for even $m$.}
	\smallskip
\noindent
For even $m$ the situation is somewhat more complicated as we shall
explain in this section.

 We take $q = 2^m$ with $m$ even, $m
\geq 4$ and we consider the code ${\cal C}_h$ with $0 < h < m /
2$. Words in ${\cal C}_h$ correspond to quadratic forms with $0 \leq w
\leq 2 h$ and $m - w \equiv 0 ({\rm mod}2)$.

In the same
way as in Proposition (3.2) we find a system of equations as condition
for $\sum_{i = 1}^{(m - w) / 2} {\rm Tr} (a_i x) {\rm Tr} (b_i x)$
to induce a word in ${\cal C}_h$. For $j = h + 1,..., (m - 2) / 2$ the
equations are of the same form as in (7) but for $j = m / 2$ we find a
term $a^{2^{m/2}} bx^{2^{m/2} + 1}$ in the expansion of ${\rm Tr} (ax)
{\rm Tr} (bx)$ which we want to neglect. To achieve that  we could
require
	$$
\sum_{i = 1}^{(m -
w)/2} a_{i}^{2^{m /2}} b_i = 0.
	$$
However, note that ${\rm Tr}
(x^{2^{m/2} + 1}) = 0$ for $x \in \F_q$ since $x^{2^{m/2}+ 1} \in
\F_{\sqrt {q}}$. So we can also disregard the term with $x^{2^{m / 2} +
1}$ if we require $\sum_{i = 1}^{(m - w)/2} a_{i}^{2^{m / 2}} b_i \in
\F_{\sqrt {q}}$ or if we require  $\sum_{i = 1}^{(m - w)/2} \bigl
(a_{i}^{2^{m / 2}} b_i + a_i b_{i}^{2^{m /2}}\bigr)= 0$ which is less
restrictive than $\sum_{i = 1}^{(m - w)/2} a_{i}^{2^{m / 2}} b_i = 0$
and of the same form as the other equations. For completeness sake we
formulate the counterpart of Proposition (3.2).
\medskip
\proclaim (5.1) Proposition.
If the elements $a_i, b_i \in \F_q$ with $1 \leq i \leq (m - w) / 2$
satisfy the system of equations
	$$
\sum_{i = 1}^{(m - w)/ 2} \bigl
(a_{i}^{2^j} b_i + a_i b_{i}^{2^j}\bigr) = 0 \eqno (18)
	$$
for $j = h +
1,..., m / 2$, then $\sum_{i = 1}^{(m - w) / 2} {\rm Tr} (a_i x) {\rm
Tr}(b_i x) $ induces a codeword in ${\cal C}_h$. \hfill\qquad $\square$
\par
\medskip
Words of minimum weight in ${\cal C}_h$ correspond to
quadratic forms with $w = 2h$ and we set $M = (m - 2h) / 2$.
Furthermore, we fix a $\F_2$-independent subset $\{a_i: 1 \leq i \leq M
\} \subset \F_q$, and  if we wish we may assume that $a_1 = 1$. The
situation is complicated by the fact that (18) has $M$ equations in $M$
unknowns $b_i$ hence we cannot guarantee the existence of $q$ solutions
$(b_1,..., b_M)$ as in the odd case. By successively applying linear
transformations to (18) we finally arrive at an equation of the form $S
(x) = 0$ with $S$ a  $2$-linearized polynomial. The actual number of
solutions of (18) depends on the choice of the subset $\{a_i: 1 \leq i
\leq M\}$.

To find generalized Hamming
weights of ${\cal C}_h$ and curves with many points we employ the same
strategy as indicated in Section 3. The reader can easily convince
himself that for even $m$ we have the analogues of  (3.5),
(3.6), (3.7) and (3.8). We shall illustrate the strategy by
working out some examples
\smallskip
\noindent I. {\sl The case } $h = (m - 2) / 2$.
\smallskip The system (18) consists of only one equation:
$$
b_{1}^{2^{m/ 2}} + b_1 = 0
$$
which has $2^{m/ 2}$ solutions or ${\rm dim}_{\F_2} (S) = m/2$. We find
by Theorem (3.6):
\medskip
\proclaim (5.2) Proposition. i). The generalized Hamming weights
of ${\cal C}_{(m - 2)/2}$ of length $2^m - 1$ satisfy $d_r ({\cal C}_{(m
- 2)/2}) = (2^r - 1) (2^{m - 1} - 2^{m - 2} )/ 2^{r - 1}$ for $1 \leq r
\leq (m - 2) / 2$.
\hfill\break
ii) For $1 \leq r \leq (m - 2) / 2$ there
exist curves $C_r$ defined over $\F_{2^m}$ of genus $g ({\cal C}_r) =
(2^r - 1) 2^{(m - 4)/ 2}$ and $\# C_r (\F_q) = 2^m + 1 +(2^r - 1) 2^{m -
1}$; these curves attain the Hasse-Weil upper bound. \hfill\qquad
$\square$
\par
\medskip
\noindent
{\bf (5.3) Example.} For $m = 6$ we have $h = 2$
and for ${\cal C}_2 = BCH (3)^{\bot}$ of length $63$ holds $d_2 (BCH
(3)^{\bot}) = {3 \over 2} d_1 (BCH(3)^{\bot}) = 24$. Furthermore there
exists curves ${\cal C}_r$ defined over $\F_{64}$ with
\medskip
\settabs 5 \columns \+&$r$ &$g ({\cal C}_r)$ &$\# C_r (\F_{64})$ &{\rm
Wirtz}\cr \smallskip
\+&$1$ &$2$ &$97$ &$97$ \cr
\smallskip
\+&$2$ &$6$ &$161$ &$155-161$\cr
\medskip
For $m = 8$ we have $h = 3$ and from Proposition (5.2.ii) we obtain
curves $C_r$ defined over $\F_{256}$ with
	\medskip
\settabs  5 \columns
\+&$r$ &$g ({\cal C}_r)$ &$\# C_r(\F_{256})$ &{\rm Upper~bound~ (1)}\cr
\smallskip
\+&$1$ &$4$ &$385$ &$385$\cr
\smallskip
\+&$2$ &$12$ &$641$ &$641$ \cr
\smallskip
\+&$3$ &$28$ &$1153$ &$1153$.\cr
\medskip
\noindent
II. {\sl The case } $h = (m - 4) / 2$.
\smallskip
We take $m \geq  6$. The system of equations is
$$
\eqalign{
b_{1}^{2^{m/2-1}} + b_1 &= a_2 b_{2}^{2^{m / 2 - 1}} + a_{2}^{2^{m /2 -
1}} b_2,\cr
b_{1}^{2^{m /2}} + b_1 &= a_2 b_{2}^{2^{m / 2}} +
a_{2}^{2^{m/2}} b_2.\cr
}
$$
If we employ the transformation $b_1 = b_{1}^\prime + \root {2^{m / 2}}
\of {a_2} b_2$ we get
$$
\eqalignno {
&b_{1}^{\prime^{2^{m/2 - 1} } } + b_{1}^\prime = \bigl (a_{2} + \sqrt
{a_2} \bigr) b_{2}^{2^{m / 2 - 1}} + \bigl (a_{2}^{2^{m / 2 - 1}} +
\root {2^{m /2}} \of {a_2}\bigr) b_2, \cr &b_{1}^{\prime^{2^{m/2}}} +
b_{1}^\prime = \bigl (a_{2}^{2^{m/2}} + \root {2^{m / 2}} \of {a_2}
\bigr) b_2 = 0. & (19) \cr}
$$
The second equation in (19) implies $b_{1}^\prime \in \F_{\sqrt {q}}$.

Now we choose $a_2 \in \F_4 - \F_2$ then the first equation in (19)
becomes
$$
b_{1}^{\prime^{2^{m /2 - 1}}} + b_{1}^\prime = b_{2}^{2^{m / 2
- 1}} + b_2.
$$
This implies $b_2 = b_{1}^\prime + \F_2$ if g.c.d. $(m , (m - 2) / 2) =
1$ or $b_2 = b_{1}^{\prime} + \F_4$ if g.c.d. $(m, (m - 2)/ 2) = 2$. We
have proved:
\medskip
\proclaim (5.4) Proposition. If we take $\{1, a_2\}$ with $a_2 \in \F_4 - \F_2$
then the system of equations (19)
has $2 \sqrt {q}$ solutions $(b_1,b_2)$ for $m \equiv 0 ({\rm mod}~4) $
and $ 4 \sqrt {q}$ solutions for $m \equiv 2 ({\rm mod}~4)$.
\hfill\qquad $\square$ \par
\medskip
We shall leave the formulation of the analogue of
(5.2) to the reader. It implies that  there exist
curves over $\F_{256}$ of genus $2$ with $321$ points and of genus $6$
with $449$ points; both values are maximal. Moreover $d_2 (BCH
(3)^{\bot}) = {3 \over 2} d_1 (BCH (3)^{\bot}) = 144$ for $BCH
(3)^{\bot}$ of length $255$.
\medskip
Finally, it is also possible to
obtain good curves by applying the method to a combination of the cases.
We shall illustrate this with an example.
\smallskip
\noindent
{\bf (5.5) Example.} Take $m = 6$ and let $\rho \in \F_4 - \F_2$. As we
saw in  case I the word ${\rm Tr} (x) {\rm Tr} (bx) \in
{\cal C}_2$ if $b^8 + b = 0$, i.e.  $b \in \F_8$. On the other hand case
II implies that ${\rm Tr} (x)  {\rm Tr}(bx) + {\rm Tr} (\rho x) {\rm Tr}
(b_2 x) \in {\cal C}_1$ if $$
\eqalignno {
&b_{1}^4 + b_1 = \rho b_{2}^4 + \rho b_2, \cr
&b_{1}^8 + b_1 = \rho b_{2}^8 + \rho^2 b_2. & (20)\cr}
$$
One checks that the solutions $(b_1, b_2)$ of (20) can be written as
$(b_1, b_2) = s (\rho, 1) + t (1, \rho)$ with $s \in \F_8, t \in
\F_4$. Now choose $b \in \F_8 - \F_2$ and $(b_1, b_2) = b(\rho, 1)$ then
$\{1, \rho, b \rho, b\}$ has $\F_2$-rank $4$. The corresponding words
are
$$
c_1 = {\rm Tr} (x) {\rm Tr} (bx) = {\rm Tr} (bx^2 + (b^2 + b) x^3 + (b^4 + b)
x^5)
\quad {\rm of~weight}\quad 16,
$$
$$c_2 = {\rm Tr} (x) {\rm Tr} (b \rho x) + {\rm Tr} (\rho x) {\rm Tr}
(bx) = {\rm Tr} ((b^2 + b) x^3)\quad {\rm of~~~weight}\quad 24.
$$

Then $c_1 + c_2 = {\rm Tr} (x) {\rm Tr} (b \rho^2 x) + {\rm Tr}(\rho x)
{\rm Tr}(bx) = {\rm Tr} (bx^2 + (b^4 +b) x^5) \in {\cal C}_2$ has weight $24$
as
follows from Proposition (2.2) since $\{1, \rho, b \rho^2, b\}$ also has
$\F_2$-rank $4$. The curves have genus $2, 1,2$ respectively and trace
of Frobenius $-32, -16, -16$ respectively. So we find a fibre product
curve over $\F_{64}$ of genus $5$ with $65+64=129$ rational points
(Wirtz: 101-145).
\smallskip
\noindent
{\bf (5.6) Remark.} For $b \in \F_8 -\F_2$ the quadratic form ${\rm
Tr}(bx^9)$ vanishes on $\F_{64}$. This implies that the curve
$y^2+y=bx^9$ with genus $4$ has trace of Frobenius $-64$. The fibre
product of this curve and
$$
y^2+y= bx^9+(b^4+b)x^5 + (b^2+b)x^3 + bx^2
$$
has genus $g=10$ and  $193$ points over $\F_{64}$ (Wirtz:
$139-225$)
\bigskip \noindent  {\bf References.}
\smallskip
\item{[A]} Arf. C.: Untersuchungen \"uber quadratische Formen in
K\"orpern der Charakteristik 2, {\sl  J. reine angew. Math.} {\bf 183}
(1941), pp. 148-167.
\item {[G-V1]} van der Geer, G. \& van der Vlugt,
M.: Reed-Muller codes and supersingular curves I, {\sl Compositio Math.}
{\bf 84} (1992), pp. 333-367.
\item {[G-V2]} van der Geer, G. \& van der
Vlugt, M.: Fibre products of Artin-Schreier curves and generalized
Hamming weights of codes, {\sl  Report  University of Amsterdam} 93-05,
1993. To appear in J. Comb. Th. A.
\item {[G-V3]} van der Geer, G. \& van der Vlugt, M.: Curves over
finite fields of characteristic 2 with many rational points, {\sl
Comptes Rendus Acad. Sci. Paris, S\'erie I}, {\bf 317} (1993), pp.
593-597.
\item {[G-V4]} van der Geer, G. \& van der Vlugt, M.:
Generalized Hamming weights of codes and curves over finite fields with
many points, {\sl Report University of Amsterdam} 94-02, 1994. To appear
in Israel Math. Conf. Proc.
\item {[L-N]} Lidl, R. \& Niederreiter,
H.: Finite fields, {\sl  Encycl. Math. Appl.}, Vol. {\bf 20}, Addison
Wesley, Reading, Mass., 1983.
\item {[S]} Serre, J-P.:  Sur le nombre
des points rationnels d'une courbe alg\'ebrique sur un corps fini,
{\sl Comptes Rendus Acad. Sci. Paris, S\'erie I}, {\bf 296} (1983), pp.
397-402 (= Oeuvres 128).
\item {[W]} Wirtz, M.: Konstruktion und
Tabellen linearer Codes, Westf\"alische Wilhelms-\hfill\break
Universit\"at M\"unster, 1991.
\bigskip
\noindent
December 1, 1994
\bigskip
\settabs3 \columns
\+G. van der Geer  &&M. van der Vlugt\cr
\+Faculteit
Wiskunde en Informatica &&Mathematisch Instituut\cr
\+Universiteit van
Amsterdam &&Rijksuniversiteit te Leiden \cr
\+Plantage Muidergracht 24&&Niels Bohrweg 1 \cr
\+1018 TV Amsterdam
&&2300 RA Leiden \cr
\+The Netherlands &&The Netherlands \cr
\bye